\documentclass[ECP]{ejpecp} 

\usepackage[round]{natbib}

\AUTHORS{%
  James~Y.~Zhao\footnote{Stanford University, USA.
    \EMAIL{jyzhao@math.stanford.edu}}
}

\SHORTTITLE{Asymptotically Ewens Measures}

\TITLE{Universality of Asymptotically\\ Ewens Measures on Partitions\thanks{Supported
    by a Stanford Graduate Fellowship.}} 


\AMSSUBJ{60F05; 60C05} 

\SUBMITTED{30 November 2012} 
\ACCEPTED{17 April 2012} 

\ARXIVID{1105.1010}  

\VOLUME{17}
\YEAR{2012}
\PAPERNUM{16}
\DOI{v17-1956}

\ABSTRACT{We give a criterion for functionals of partitions to converge to a universal limit under a class of measures that ``behaves like'' the Ewens measure. Various limit theorems for the Ewens measure, most notably the Poisson-Dirichlet limit for the longest parts, the functional central limit theorem for the number of parts, and the Erd\H os-Tur\'an limit for the product of parts, extend to these asymptotically Ewens measures as easy corollaries. Our major contributions are: (1) extending the classes of measures for which these limit theorems hold; (2) characterising universality by an intuitive and easily-checked criterion; and (3) providing a new and much shorter proof of the limit theorems by taking advantage of the Feller coupling.}

\renewcommand\P{\mathbb P}
\newcommand\E{\mathbb E}

\newcommand\R{\mathbb R}

\newcommand\N{\mathbb N}
\newcommand\I{\mathbb 1}

\newcommand\cD{\mathcal D}
\newcommand\cL{\mathcal L}
\newcommand\cP{\mathcal P}

\newcommand\cX{\mathcal X}

\newcommand\cites[1]{\citeauthor{#1}'s (\citeyear{#1})}

\begin{document}



\section{Introduction}
Let $\cP_n$ be the partitions of $n\in\N$, which we represent by $\alpha=(\alpha_1,\ldots,\alpha_n)\in\cP_n$, where $\alpha_i$ is the number of parts of size $i$, so that $\alpha_1+2\alpha_2+\cdots+n\alpha_n=n$. Let 
\begin{equation}
\P^\theta_n(\alpha)=\frac{n!}{\theta(\theta+1)\cdots(\theta+n-1)}\prod_{i=1}^n\frac{\theta^{\alpha_i}}{i^{\alpha_i}\alpha_i!}
\end{equation}
be the \citet{ewens} measure on $\cP_n$ with parameter $\theta>0$. The Ewens measure should be interpreted as a random permutation weighted by the number of cycles, under which the $\alpha_i$ are asymptotically independent Poisson with parameter $\theta/i$. Any other probability measure $\P_n$ on $\cP_n$ can be written as a perturbation of the Ewens measure with parameter $\theta>0$, that is,
\begin{equation}
\P_n(\alpha)=\P^{\theta,\eta}_n(\alpha)=\frac{\eta(\alpha)}{Z_n^{\theta,\eta}}\prod_{i=1}^n\frac{\theta^{\alpha_i}}{i^{\alpha_i}\alpha_i!},
\end{equation}
where $\eta:\cP_n\rightarrow\R^+$ is (any multiple of) the Radon-Nikodym derivative of $\P_n$ with respect to the Ewens measure $\P^\theta_n$, and $Z^{\theta,\eta}_n$ is a normalising constant. Since $n$, $\theta$ and $\eta$ uniquely determine the measure, it will be convenient to use the notation $\P^{\theta,\eta}_n$ to specify it.

We will study partitions by their preimages under the \cite{feller} coupling, which we define as the map $\{0,1\}^\N\rightarrow\cP_n$ which sends a binary sequence $(\xi_1,\xi_2,\ldots)\in\{0,1\}^\N$ to a partition $(\alpha_1(n),\ldots,\alpha_n(n))\in\cP_n$, where
\begin{equation}
\alpha_i(n)=\sum_{j=1}^{n-i}\xi_j(1-\xi_{j+1})\cdots(1-\xi_{j+i-1})\xi_{i+j}+\xi_{n-i+1}(1-\xi_{n-i+2})\cdots(1-\xi_n)
\end{equation}
is the number of gaps between 1s of length $i-1$ in the string $(\xi_1,\xi_2,\ldots,\xi_n,1)$. It was shown by \cite{ABT1992} that if $\P^\theta_F$ is a measure on $\{0,1\}^\N$ under which the $\xi_i$ are independent Bernoulli variables with probability of success $\theta/(\theta+i-1)$, then the measure induced on $\cP_n$ by the Feller coupling is the Ewens measure $\P^\theta_n$. In particular, we can embed $\P^\theta_n$ within $\P^\theta_F$ so that $\alpha_i=\alpha_i(n)$, simultaneously for all $n\in\N$.

For any $m,n\in\N$, $(\alpha_1(n),\ldots,\alpha_m(n))$ is a partition of some integer $k\le n$. Then, for any sequence of measures $\P^{\theta,\eta}_n$ on $\cP_n$, $n\in\N$, define the \textit{weights}
\begin{equation}
\eta_{n,m}=\eta(\alpha_1(n),\ldots,\alpha_m(n)),
\end{equation}
where $\eta$ is the (unscaled) Radon-Nikodym derivative of $\P^{\theta,\eta}_k$ with respect to $\P^\theta_k$, for the appropriate choice of $k=\alpha_1(n)+2\alpha_2(n)+\cdots+m\alpha_m(n)$. If we interpret the Radon-Nikodym derivative as weights, then $\eta_{n,m}$ should be interpreted as the weight of a partition of $n$ considering only parts of size at most $m$. Note that since the weights are defined in terms of the $\alpha_i(n)$, they are random variables under the Feller measure $\P^\theta_F$.

We call a sequence of measures $\P^{\theta,\eta}_n$ \textit{asymptotically Ewens} when:
\begin{enumerate}
\item The limits $\displaystyle\lim_{n\rightarrow\infty}\eta_{n,n}$ and $\displaystyle\lim_{m\rightarrow\infty}\lim_{n\rightarrow\infty}\eta_{n,m}$ exist and agree in $L^1(\P^\theta_F)$; and
\item The common limit $\eta_\infty$ satisfies $0<\E_F^\theta\big[\eta_\infty\big]<\infty$.
\end{enumerate}

These measures are important in a variety of applications and include many extensively studied measures as special cases. Essentially, they generalise the \textit{logarithmic combinatorial structures} of \citet{ABT2000} by removing the conditioning relation. We will prove this claim in Section \ref{sec:lcs}, while Section \ref{sec:background} contains background on these measures and the associated limit theorems.

We conclude the introduction by stating our main theorem, the proof of which is presented in Section \ref{sec:theorem}, and showing how the Ewens measure limit theorems extend to asymptotically Ewens measures as easy corollaries.

\begin{theorem}
\label{thm:universality}
Suppose $X_n:\cP_n\rightarrow(\cX,||\cdot||)$ is a sequence of deterministic functions on partitions of $n$ (and therefore they are functions on $\{0,1\}^n$ by the Feller coupling), taking values in some normed space, such that for any fixed $d\in\N$,
\begin{equation}
\label{eqn:mainassumption}
\lim_{n\rightarrow\infty}\max_{\xi_1,\ldots,\xi_n}\Big|\Big|X_n(\xi_1,\ldots,\xi_n)-X_n(1,\ldots,1,\xi_{d+1},\ldots,\xi_n)\Big|\Big|=0.
\end{equation}
If $X_n\stackrel d\longrightarrow X$ under the Ewens measure with parameter $\theta$, then $X_n\stackrel d\longrightarrow X$ under any asymptotically Ewens measure with parameter $\theta$.
\end{theorem}
The criterion (\ref{eqn:mainassumption}) has a very intuitive interpretation: if the functional does not depend on the first finitely many Feller variables, then it has a universal limit for any asymptotically Ewens measure, depending only on the parameter $\theta$. 

\begin{corollary}[Poisson-Dirichlet]
\label{cor:PD}
Let $L_{n,k}$ be the length of the $k$th longest part of a partition of $n$, and let $L_n=(L_{n,1},L_{n,2},\ldots)$. Under any asymptotically Ewens measure with parameter $\theta$, $L_n/n$ converges in distribution in $L^1(\N)$ to $PD(\theta)$, the Poisson-Dirchlet measure with parameter $\theta$ \citep{kingman1975,kingman1977}.
\end{corollary}
\begin{proof}
Consider strings $(1,\ldots,1,\xi_{i+1},\ldots,\xi_n)$, for $i=d,d-1,\ldots,0$. For each decrement of $i$, either the partition is unchanged, or a 1-part is deleted and some other part length increases by 1, which does not change the order of parts and thus changes $L_n$ by at most 2 in the $L^1(\N)$ norm. Since there are $d$ decrements from $i=d$ to $i=0$,
\begin{equation}
\max_{\xi_1,\ldots,\xi_n}\Big|\Big|L_n(\xi_1,\ldots,\xi_n)-L_n(1,\ldots,1,\xi_{d+1},\ldots,\xi_n)\Big|\Big|_1\le2d.
\end{equation}
Hence, $X_n=L_n/n$ satisfies the conditions of Theorem \ref{thm:universality} applied to the normed space $L^1(\N)$, so the result follows from the Poisson-Dirichlet limit for the Ewens measure \citep{kingman1975,kingman1977,watterson}.
\end{proof}

\begin{corollary}[CLT]
\label{cor:FCLT}
Let $\nu_{n,t}=\alpha_1+\cdots+\alpha_{\lfloor n^t\rfloor}$, $0\le t\le1$, be the number of parts of size at most $n^t$ in a partition of $n$. Under any asymptotically Ewens measure with parameter $\theta$, $(\nu_{n,t}-\theta t\log n)/\sqrt{\theta \log n}$ converges in distribution in $\cD[0,1]$, the Skorohod space of right-continuous left-limit functions on $[0,1]$, to the standard Brownian motion $W_t$.
\end{corollary}
\begin{proof}
Consider strings $(1,\ldots,1,\xi_{i+1},\ldots,\xi_n)$, for $i=d,d-1,\ldots,0$. For each decrement of $i$, either the partition is unchanged, or a 1-part is deleted and some other part length increases by 1, which changes $\nu_{n,t}$ by at most 1. Since there are $d$ decrements,
\begin{equation}
\max_{\xi_1,\ldots,\xi_n}\Big|\nu_{n,t}(\xi_1,\ldots,\xi_n)-\nu_{n,t}(1,\ldots,1,\xi_{d+1},\ldots,\xi_n)\Big|\le d.
\end{equation}
Letting $\nu_n\in\cD[0,1]$ be the sample path of $\nu_{n,t}$ for $0\le t\le 1$,
\begin{equation}
\max_{\xi_1,\ldots,\xi_n}\Big|\Big|\nu_n(\xi_1,\ldots,\xi_n)-\nu_n(1,\ldots,1,\xi_{d+1},\ldots,\xi_n)\Big|\Big|_\infty\le d.
\end{equation}
Hence, $X_n=(\nu_n-\theta t\log n)/\sqrt{\theta \log n}$ satisfies the conditions of Theorem \ref{thm:universality} applied to the normed space $\cD[0,1]$, so the result follows from the functional central limit theorem for the Ewens measure \citep{delaurentispittel,hansen1990} .
\end{proof}

\begin{corollary}[Erd\H os-Tur\'an]
\label{cor:ET}
Let $O_{n,t}=lcm\big\{i\le n^t:\alpha_i>0\big\}$ be the least common multiple (or product) of the parts of size at most $n^t$ in a partition of $n$. Under any asymptotically Ewens measure with parameter $\theta$, $\big(\!\log O_{n,t}-\theta t^2(\log n)^2/2\big)/\sqrt{\theta(\log n)^3/3}$ converges in distribution in $\cD[0,1]$ to $W_{t^3}$, where $W_t$ is the standard Brownian motion.
\end{corollary}
\begin{proof}
The partitions given by $(\xi_1,\ldots,\xi_n)$ and $(1,\ldots,1,\xi_{d+1},\ldots,\xi_n)$ differ in at most $d$ parts. Since deleting a part of size $\ell$, adding a part of size $m$, or replacing $\ell$ by $m$ changes the logarithm of the least common multiple (or product) by at most $\log\max(\ell,m)$,
\begin{equation}
\max_{\xi_1,\ldots,\xi_n}\Big|\log O_{n,t}(\xi_1,\ldots,\xi_n)-\log O_{n,t}(1,\ldots,1,\xi_{d+1},\ldots,\xi_n)\Big|\le d\log n^t.
\end{equation}
Then, letting $O_n\in\cD[0,1]$ be the sample path of $O_{n,t}$ for $0\le t\le1$,
\begin{equation}
\max_{\xi_1,\ldots,\xi_n}\Big|\Big|\log O_n(\xi_1,\ldots,\xi_n)-\log O_n(1,\ldots,1,\xi_{d+1},\ldots,\xi_n)\Big|\Big|_\infty\le d\log n.
\end{equation}
Hence, $X_n=\big(\log O_n-\theta t^2(\log n)^2/2\big)/\sqrt{\theta(\log n)^3/3}$ satisfies the conditions of Theorem \ref{thm:universality}, so the result follows from the functional Erd\H os-Tur\'an limit for the Ewens measure \citep{erdosturan1,erdosturan3,ABT1994}. The same result holds for the product of parts.
\end{proof}

\section{Background}
\label{sec:background}
\subsection{Measures on Partitions}
Measures on partitions arise naturally from combinatorial objects which consist of components of various sizes. For example, cycles of a random permutation, irreducible factors of a random polynomial or Jordan blocks of a random matrix are all described by measures on partitions when one cares about only the sizes of those components. More such examples are given by \citet{ABT1997,ABTbook}.

The most basic example is a uniformly random permutation, which corresponds to our $\P^{\theta,\eta}_n$ when $\theta=1$ and $\eta=1$ identically, and has been the subject of extensive study since the 19th century. The generalisation to $\theta>0$, still with $\eta=1$ identically, was introduced by \citet{ewens} to model propagation of alleles in population genetics, and represents a random permutation weighted by the number of cycles, or perhaps more intuitively, a permutation formed in a Markov process where cycles are added at a rate $\theta$ \citep{hoppe}.

The further generalisation to weights $\eta(\alpha)=\zeta_1(\alpha_1)\zeta_2(\alpha_2)\cdots\zeta_n(\alpha_n)$, with certain limiting conditions on the $\zeta_i$, was introduced by \citet{ABT2000}. Their \textit{logarithmic combinatorial structures} generalise the \textit{decomposable combinatorial structures} of \citet{flajoletsoria}, which are measures on partitions induced by the uniform measure on families of combinatorial objects determined by sizes of components. We will prove in Section \ref{sec:lcs} that logarithmic combinatorial structures are exactly asymptotically Ewens measures with weight function in the form above.

An important subclass of logarithmic combinatorial structures are measures with weights $\eta(\alpha)=\zeta_1^{\alpha_1}\zeta_2^{\alpha_2}\cdots\zeta_n^{\alpha_n}$, where $\zeta_i$ are constants with $\sum_i|\zeta_i-1|/i<\infty$. This is the asymptotically Ewens case of the \textit{multiplicatively weighted measures} of \cites{betz2009} combinatorial model for Bose-Einstein condensation. In this model, $n$ points are determined by both their positions and a permutation that describes their trajectories, with the energy function (Hamiltonian) of the system depending particularly on the presence of long cycles. \citet{betz2010} analysed the asymptotic behaviour of the cycle lengths and the normalising constant (partition function) in the asymptotically Ewens case, as well as two other cases where the weights diverge. \citet{ercolani} continued this work, extending the cycle length analysis to a much more comprehensive set of parameter regions.

The further specialisation $\zeta_i=(1-q^i)/(1-t^i)$, for parameters $0<q<1$ and $0<t<1$, are the \textit{Macdonald polynomial measures} of \cites{diaconisram} probabilistic interpretation of Macdonald polynomials. This measure is the stationary distribution of a Markov chain on partitions which corresponds to the Macdonald operator on symmetric polynomials.

Some examples of asympototically Ewens measures which are not logarithmic combinatorial structures include many restricted permutations, such as permutations with more even cycles than odd cycles, permutations whose squares have fixed points, or permutations with an even number of cycles; for any underlying measure that is a logarithmic combinatorial structure, these restrictions are asymptotically Ewens.

There are some examples of measures, such as the $a$-riffle shuffle measures of \citet{diaconismcgrathpitman}, and the restricted permutations studied by \citet{lugo}, which are not asymptotically Ewens by our current definition, but behave similarly in the sense that they follow the Poisson-Dirichlet limit, as discussed in more detail in Section \ref{sec:PD}. Generalising our measures to include these examples would be an interesting direction for future work.

Finally, there are many measures which are not asymptotically Ewens in any sense, such as the uniform measure on partitions, \cites{pitman} two-parameter family of measures (although this family includes the Ewens measure as a special case), and the induced measure on partitions from various measures on permutations such as the Plancherel measure and its generalisation, the Schur measures of \citet{okounkov}.

\subsection{The Poisson-Dirichlet Limit}
\label{sec:PD}
The parameter $\theta$ in the Ewens measure corresponds to a rate of formation of new parts \citep{hoppe}; indeed, $\theta$ is the global rate of mutation in \citeauthor{ewens}' (\citeyear{ewens}) original genetic model. This insight extends to the asymptotically Ewens case, where the rate of formation of new parts, appropriately scaled, converges to the parameter $\theta$. This is the intuitive reason why we expect the limit theorems to be universal: with new parts being added at the same rate, the relative sizes of parts should behave similarly.

The key notion that captures this behaviour is the Poisson-Dirichlet limit: the largest parts, normalised by $\frac1n$, converge in distribution on $L^1(\N)$ to a limit known as the Poisson-Dirichlet measure with parameter $\theta$. This was first studied by \citet{kingman1975}, who described it via the Dirichlet distribution on $L^1(\N)$, and \citet{watterson}, who found an explicit density.

Historically, \citet{golomb} was the first to calculate the expected value of the longest cycle of a uniformly random permutation, \citet{shepplloyd} found the distributions of the $k$th longest cycles, and \citet{kingman1975,kingman1977} and \citet{watterson} found the joint distribution of longest cycles under the Ewens measure. \citet{hansen1994} proved the Poisson-Dirichlet limit for decomposable combinatorial structures, while the extension to logarithmic combinatorial structures was made by \citet{ABT1999}.

We further generalise the Poisson-Dirichlet limit to asymptotically Ewens measures. However, there are still many other measures which satisfy the Poisson-Dirichlet limit, such as the largest prime factors of a random integer \citep{knuth}, the $a$-riffle shuffle measures of \citet{diaconismcgrathpitman}, and the restricted permutations studied by \citet{lugo}. We expect there to be a fundamental reason why we observe the same limit in these cases, although what that reason should be is currently beyond our grasp.
 
The Poisson-Dirichlet distribution has a two-parameter generalisation \citep{pitmanyor}, which is the limit of the ordered parts of \cites{pitman} two-parameter family of measures on partitions. Since Pitman's measures are a direct generalisation of the Ewens measure, it seems plausible that our result could be extended in this direction.

\subsection{Other Limit Theorems}
\label{sec:otherlimits}
It is classical that the number of $i$-cycles in a uniformly random permutation are asymptotically independent Poisson with parameter $1/i$. The usual proof is by generating functions, which also works for multiplicative weights $\eta(\alpha)=\prod_i\zeta_i^{\alpha_i}$, where the parameter becomes $\theta\zeta_i/i$. Such a proof is given by \citet{betz2010}; see also the book of \citet{ABTbook} for a thorough treatment of generating function techniques in this setting. For logarithmic combinatorial structures, \citet{ABT2000} prove that the number of parts of size $i$ are asymptotically independent, although this result is in some sense one of the defining assumptions of logarithmic combinatorial structures.

The total number of parts was first studied by \citet{goncharov1942}, who found a central limit theorem for the number of cycles in a uniformly random permutation. The functional central limit theorem as seen in Corollary \ref{cor:FCLT} was first proved by \citet{delaurentispittel} for the uniform permutation case, and extended to the Ewens measure by \citet{hansen1990}. \citet{flajoletsoria} proved a central limit theorem for decomposable combinatorial structures, and the two theorems were unified by \citet{ABT2000}, who proved a functional central limit theorem for logarithmic combinatorial structures.

The asymptotic moments of the smallest parts were derived by \citet{shepplloyd}. They have not been the subject of extensive study; some facts which are known about them are listed in the book of \citet{ABTbook}. The shortest parts depend heavily on the first few Feller variables, and do not fall under the scope of our universality theorem.

It is also possible to canonically order the parts under the Ewens measure by the order in which they appear in the \textit{Chinese restaurant coupling}. This limit is called the Griffiths-Engen-McCloskey measure \citep{griffiths}, and behaves similarly to the longest parts; in fact its order statistics exactly follow the Poisson-Dirichlet measure. This limit theorem does not generalise to asymptotically Ewens measures due to the a lack of a canonical order.

The lowest common multiple of parts is a statistic of interest when the partition is induced by a permutation, as it is the group order of the permutation. \citet{erdosturan1,erdosturan3} found a central limit theorem for the logarithm of the lowest commmon multiple in the uniform permutation case. Their proof, and all subsequent proofs, worked via the product of parts, in particular showing that the product satisfies the same central limit theorem. The generalisation to the Ewens measure was proved by \citet{ABT1994}, who also proved the functional form in Corollary \ref{cor:ET}, while the extension to logarithmic combinatorial structures was made by \citet{ABT2000}.

\subsection{Universality}
For logarithmic combinatorial structures, \citet{ABT2000} prove
\begin{equation}
\label{eqn:ABTtheorem}
\lim_{n\rightarrow\infty}\big|\big|\cL^{\theta,\eta}_n(\alpha_{d_n},\ldots,\alpha_n)-\cL^\theta_n(\alpha_{d_n},\ldots,\alpha_n)\big|\big|_{TV}=0,
\end{equation}
where $\cL^{\theta,\eta}_n$ and $\cL^\theta_n$ are the laws of $(\alpha_{d_n},\ldots,\alpha_n)$ under $\P^{\theta,\eta}_n$ and $\P^\theta_n$ respectively, and $d_n$ is any sequence satisfying $d_n\rightarrow\infty$ and $d_n/n\rightarrow0$. They also proved the other limit theorems above, but most of their proofs ran in parallel to (\ref{eqn:ABTtheorem}) instead of directly invoking it.

This left open the question of a simple criterion to determine whether a functional is universal, as well as the question of whether asymptotic independence of the $\alpha_i$ is a necessary condition. Our model answers both of these questions, removing the requirement for the $\alpha_i$ to be asymptotically independent, and also giving an easily-checked criterion for universality.

\section{Proof of Main Theorem}
\label{sec:theorem}
Suppose $\P^{\theta,\eta}_n$ is an asymptotically Ewens measure with parameter $\theta>0$ and $X_n$ satisfies the conditions of Theorem \ref{thm:universality}. It suffices \citep{billingsley} to prove that for any bounded, uniformly continuous function $f:(\cX,||\cdot||)\rightarrow(\R,|\cdot|)$, we have $\E^{\theta,\eta}_n\big[f(X_n)\big]\sim\E^\theta_n\big[f(X_n)\big]$ as $n\rightarrow\infty$, where $a_n\sim b_n$ means $a_n/b_n\rightarrow1$.

Since $\eta_{n,n}$ is (a multiple of) the Radon-Nikodym derivative,
\begin{equation}
\label{eqn:mainthm1}
\E_n^{\theta,\eta}\big[f(X_n)\big]
=\frac{Z_n^\theta}{Z_n^{\theta,\eta}}\E_F^\theta\big[\eta_{n,n}f(X_n)\big].
\end{equation}
Since $f$ is bounded and $\eta$ is asymptotically Ewens, $(\eta_{n,n}-\eta_{n,m})f(X_n)\rightarrow0$ in $L^1(\P^\theta_F)$ as $n\rightarrow\infty$ then $m\rightarrow\infty$, hence
\begin{equation}
\label{eqn:mainthm2}
\lim_{n\rightarrow\infty}\E _F^\theta\big[\eta_{n,n}f(X_n)\big]
=\lim_{m\rightarrow\infty}\lim_{n\rightarrow\infty}\E_F^\theta\big[\eta_{n,m}f(X_n)\big].
\end{equation}
For any integer $d<n$, the event $\{\alpha_i(n)\ne\alpha_i(d)\}$ can be decomposed into a union of two events: either $(\xi_{d-i+1},\ldots,\xi_{d+1})=(1,0,\ldots,0)$, with probability
\begin{equation}
\E^\theta_F\big[\xi_{d-i+1}(1-\xi_{d-i+2})\cdots(1-\xi_{d+1})\big]\le\E_F^\theta\big[\xi_{d-i+1}\big]=\tfrac\theta{d-i+\theta},
\end{equation}
or $(\xi_{d-i+2},\ldots,\xi_n,1)$ contains a substring $(1,0,\ldots,0,1)$ with $i-1$ zeroes, with probability bounded above by
\begin{equation}
\sum_{k=2}^{n-d}\E^\theta_F\big[\xi_{d-i+k}\xi_{d+k}\big]+\E^\theta_F\big[\xi_{n-i+1}\big]
=\sum_{k=2}^{n-d}\tfrac\theta{d-i+k-1+\theta}\tfrac\theta{d+k-1+\theta}+\tfrac\theta{n-i+\theta}.
\end{equation}
Using telescoping series, we have $\P^\theta_F[\alpha_i(n)\ne\alpha_i(d)]\le(2\theta+\theta^2)/(d-i+\theta)$. Hence, the event $E\equiv\{\exists\,i\le m:\alpha_i(n)\ne\alpha_i(d)\}$ has $\P^\theta_F$ probability at most $(2\theta+\theta^2)m/(d-m+\theta)$, which converges to 0 as $d\rightarrow\infty$. Note that
\begin{equation}
\eta_{n,m}f(X_n)=\eta_{d,m}f(X_n)+(\eta_{n,m}-\eta_{d,m})f(X_n)\I_E.
\end{equation}
Since $\eta_{n,m}$ and $\eta_{d,m}$ converge in $L^1(\P^\theta_F)$, they are uniformly integrable. Since $\I_E\rightarrow0$ in probability in the limit $n\rightarrow\infty$ then $d\rightarrow\infty$ then $m\rightarrow\infty$, $(\eta_{n,m}-\eta_{d,m})f(X_n)\I_E\rightarrow0$ in probability and therefore in $L^1(\P^\theta_F)$, hence
\begin{equation}
\label{eqn:mainthm3}
\lim_{m\rightarrow\infty}\lim_{n\rightarrow\infty}\E_F^\theta\big[\eta_{n,m}f(X_n)\big]
=\lim_{m\rightarrow\infty}\lim_{d\rightarrow\infty}\lim_{n\rightarrow\infty}\E_F^\theta\big[\eta_{d,m}f(X_n)\big].
\end{equation}
Uniform continuity of $f$ and (\ref{eqn:mainassumption}) imply $\E^\theta_F\big[f(X_n)\big|\xi_1,\ldots,\xi_d\big]\sim\E^\theta_F\big[f(X_n)\big]$ uniformly as $n\rightarrow\infty$. Since $f$ is bounded and $\eta_{d,m}$ is uniformly integrable,
\begin{align}
\label{eqn:mainthm4}
\lim_{n\rightarrow\infty}\E^\theta_F\big[\eta_{d,m}f(X_n)\big]
&=\lim_{n\rightarrow\infty}\E^\theta_F\Big[\eta_{d,m}\E^\theta_F\big[f(X_n)\big|\xi_1,\ldots,\xi_d\big]\Big]\\
&=\lim_{n\rightarrow\infty}\E^\theta_F\big[\eta_{d,m}\big]\E^\theta_F\big[f(X_n)\big].\nonumber
\end{align}
Finally, as $n\rightarrow\infty$ then $d\rightarrow\infty$ then $m\rightarrow\infty$,
\begin{equation}
\label{eqn:mainthm5}
\E^\theta_F\big[\eta_{d,m}\big]\sim\E^\theta_F\big[\eta_{n,n}\big]
=\sum_{\alpha\in\cP_n}\frac{\eta(\alpha)}{Z^\theta_n}\prod_{i=1}^n\frac{\theta^{\alpha_i}}{i^{\alpha_i}\alpha_i!}
=\frac{Z^{\theta,\eta}_n}{Z^\theta_n}.
\end{equation}
We have used the assumption that this quantity has a positive, finite limit. The chain of asymptotic equivalences in (\ref{eqn:mainthm1}), (\ref{eqn:mainthm2}), (\ref{eqn:mainthm3}), (\ref{eqn:mainthm4}) and (\ref{eqn:mainthm5}) completes the proof.
\hfill\qedsymbol

\section{Logarithmic Combinatorial Structures}
\label{sec:lcs}
A \textit{uniform logarithmic combinatorial structure} \citep{ABT2000} is a sequence of measures $\P_n$ on $\cP_n$, $n\in\N$, such that:
\begin{enumerate}
\item (\textit{Conditioning Relation}) There is some sequence of independent random variables $Y_1,Y_2,\ldots$ such that $\P_n(\alpha)=\P\big[\forall\,i\le n,Y_i=\alpha_i\,\big|\sum_{i\le n}iY_i=n\big]$; and
\item (\textit{Uniform Logarithmic Condition}) Each variable $Y_i$ satisfies $\big|i\P[Y_i=1]-\theta\big|\le e_i$ and $i\P[Y_i=\ell]\le e_ic_\ell$ for $\ell\ge2$, where $e_i$ and $c_\ell$ are vanishing sequences such that $e_i/i$ and $\ell c_\ell$ are summable.
\end{enumerate}

\begin{lemma}
\label{lem:lcs}
Any uniform logarithmic combinatorial structure $\P_n$ can be written as $\P^{\theta,\eta}_n$ for $\eta(\alpha)=\prod_i\zeta_i(\alpha_i)$, where $\zeta_i(0)=1$, $\big|\zeta_i(1)-1\big|\le e_i/\theta$ and $\zeta_i(\ell)\le i^{\ell-1}\ell!e_ic_\ell/\theta^\ell$ for $\ell\ge2$. As before, $e_i$ and $c_\ell$ are vanishing sequences with $e_i/i$ and $\ell c_\ell$ summable. Additionally, with $c_0=0$ and $c_1=1$, we can insist that for each $i$, $\I_{\{\ell\le1\}}+i^{\ell-1}\ell!e_ic_\ell/\theta^\ell$ is monotonic increasing in $\ell\ge0$.
\end{lemma}
\begin{proof}
Let $p_i=\P\big[Y_i=0\big]$, and let $\zeta_i(\ell)=\P\big[Y_i=\ell\big](i/\theta)^\ell\ell!/p_i$. By the conditioning relation, $\P_n=\P^{\theta,\eta}_n$, $\zeta_i(0)=1$. By the uniform logarithmic condition, $\big|p_i\zeta_i(1)-1\big|\le e_i/\theta$ and $\zeta_i(\ell)\le i^{\ell-1}\ell!e_ic_\ell/p_i\theta^\ell$ for $\ell\ge2$.

Note that $1-p_i=\sum_{\ell\ge1}\P\big[Y_i=\ell\big]\le(\theta+Ce_i)/i$, where $C=\sum_\ell c_\ell<\infty$. In particular, $p_i\rightarrow1$ as $i\rightarrow\infty$, so
\begin{equation}
\big|\zeta_i(1)-1\big|
\le(1-p_i)\zeta_i(1)+\big|p_i\zeta_i(1)-1\big|
=O(\tfrac1i)+O(e_i).
\end{equation}
Thus, $e'_i=\max\big(e_i/p_i,\theta\big|\zeta_i(1)-1\big|,1/i\big)$ is a vanishing sequence such that $e'_i/i$ is summable, and we have the required inequalities $\big|\zeta_i(1)-1\big|\le e'_i/\theta$ and $\zeta_i(\ell)\le i^{\ell-1}\ell!e'_ic_\ell/\theta^\ell$ for $\ell\ge2$. It remains to replace $c^{}_\ell$ by $c'_\ell$ satisfying the desired monotonicity condition.

Let $c'_0=0$, $c'_1=1$, and $c'_2=\max\big(c_2,\sup_i(\theta^2+\theta e'_i)/(2ie'_i)\big)$, which is finite since $e'_i\ge1/i$. For $3\le\ell\le2\theta$, let $c'_\ell=\max(c_\ell,\theta c'_{\ell-1})$, and for $\ell>2\theta$, let $c'_\ell=\max\big(c_\ell,\frac12c'_{\ell-1}\big)$. There are $\ell_0=\lfloor2\theta\rfloor<\ell_1<\ell_2<\cdots$ such that for $\ell_0\le\ell<\ell_1$, $c'_\ell=2^{\ell_0-\ell}c'_{\ell_0}$, and for $\ell_k\le\ell<\ell_{k+1}$, $c'_\ell=2^{\ell_k-\ell}c_{\ell_k}$, hence
\begin{equation}
\sum_{\ell\ge2\theta}\ell c'_\ell
\le\bigg(\sum_{j\ge0}(1+j)2^{-j}\bigg)\bigg(\ell_0c'_{\ell_0}+\sum_{k\ge1}\ell_kc_{\ell_k}\bigg)<\infty.
\end{equation}
We also have $c^{}_\ell\le c'_\ell$ and the monotonicity condition by construction.
\end{proof}

\begin{theorem}
\label{thm:lcs}
Uniform logarithmic combinatorial structures are asymptotically Ewens.
\end{theorem}
\begin{proof}
Let $\alpha_i(\infty)=\sum_{j=1}^\infty\xi_j(1-\xi_{j+1})\cdots(1-\xi_{j+i-1})\xi_{i+j}$ be the number of gaps between 1s of size $i-1$ in the string $(\xi_1,\xi_2,\ldots)$, which are independent and Poisson with parameter $\theta/i$ \citep{ABT1992}.

We will prove $\eta_{n,m}$ and $\eta_{n,n}$ are uniformly integrable and converge in probability to $\eta_\infty=\prod_i\zeta_i(\alpha_i(\infty))$, which has positive and finite expectation, by defining several intermediate weights and proving a sequence of asymptotic equivalences between them.
\begin{itemize}
\item Let $\eta_{\infty,m}=\prod_{i=1}^m\zeta_i(\alpha_i(\infty))$.
\item Let $\eta^+_{n,m}$ be weights for $\zeta^+_i(\ell)=\I_{\{\ell\le1\}}+i^{\ell-1}\ell!e_ic_\ell/\theta^\ell$, and let $\eta^-_{n,m}$ be weights for $\zeta^-(0)=1$, $\zeta^-(1)=\max(1-e_i/\theta,0)$ and $\zeta^-(\ell)=0$ for $\ell\ge2$. By Lemma \ref{lem:lcs}, we have $\eta^-_{n,m}\le\eta^{}_{n,m}\le\eta^+_{n,m}$, and $\eta^+_{n,m}$ is monotonic increasing in each $\alpha_i(n)$.
\item Let $\tilde \eta_{n,m}=\prod_{i=1}^m\zeta_i(\tilde\alpha_i(n))$ be the weight given by replacing $\alpha_i(n)$ by $\tilde\alpha_i(n)$, where $\tilde\alpha_i(n)$ is the number of gaps between 1s of length $i-1$ in the string $(\xi_1,\ldots,\xi_n,0)$.
\end{itemize}
For brevity, all limits are implicitly $n\rightarrow\infty$ with $m$ fixed, then $m\rightarrow\infty$. We also omit writing the measure explicitly; the only measure used is $\P^\theta_F$.

\vspace{6pt}\noindent\textbf{Step 1: Positivity and finiteness of limit.} Let $C=\sum_\ell c_\ell<\infty$. Since the $\alpha_i(\infty)$, $i\in\N$, are independent and Poisson with parameter $\theta/i$, we can calculate
\begin{align}
\E\big[\eta^+_\infty\big]
&=\prod_{i=1}^\infty\E\big[\zeta^+_i(\alpha_i(\infty))\big]
=\prod_{i=1}^\infty e^{-\frac\theta i}\bigg(1+\frac{\theta+Ce_i}i\bigg)\\
&\le\prod_{i=1}^\infty\bigg(1+\frac{Ce_i}i\bigg)
=\exp\sum_{i=1}^\infty\log\bigg(1+\frac{Ce_i}i\bigg).\label{eqn:etaplus}
\end{align}
We have used the inequality $1+x+y\le e^x(1+y)$. Since $e_i/i$ is summable, so is the series above, hence $\E\big[\eta^+_\infty\big]<\infty$. A similar calculation shows that $\E\big[\eta^-_\infty]>0$. 
Since $\eta^-_\infty\le\eta^{\vphantom+}_\infty\le\eta^+_\infty$, it follows that $0<\E[\eta_\infty]<\infty$.

\vspace{6pt}\noindent\textbf{Step 2: Convergence in probability.} We will show that
\begin{equation}
\eta_{n,n}\stackrel p\sim\tilde \eta_{n,n}\stackrel p\sim\tilde \eta_{n,m}\stackrel p\sim \eta_{n,m}\stackrel p\sim \eta_{\infty,m}\stackrel p\sim \eta_\infty,
\end{equation}
where $X_n\stackrel p\sim Y_n$ means $X_n/Y_n\stackrel p\rightarrow1$. The relation $\stackrel p\sim$ is clearly transitive; also observe that $\eta_{n,n}\stackrel p\sim \eta_\infty$ implies $\eta_{n,n}\stackrel p\rightarrow \eta_\infty$, since for any $\epsilon>0$, there exists some $M$ such that $\P\big[\eta_\infty>M\big]<\frac\epsilon2$, and there exists some $N$ depending on $M$ such that for all $n>N$, $\P\big[\big|\eta_{n,n}/\eta_\infty-1\big|>\frac\epsilon M\big]<\frac\epsilon2$.

Observe that for each fixed $n$, $\alpha_i(n)$ and $\tilde\alpha_i(n)$ are equal for all $i$ except one value $i^*(n)$ where $\alpha_{i^*(n)}(n)=\tilde\alpha_{i^*(n)}(n)+1$, and given $i^*(n)=i$, $\tilde\alpha_i(n)$ has the same law as an independent copy of $\alpha_i(n-i)$. Furthermore, $\alpha_i(n-i)\ne0$ implies either $\alpha_i(\infty)\ne0$, with probability $1-e^{-\theta/i}$, or $i^*(n-i)=i$, with probability $\frac1{n-i}$. Hence,
\begin{equation}
\P\big[\tilde\alpha_{i^*(n)}(n)\ne0\big]
=\frac1n\sum_{i=1}^n\P\big[\alpha_i(n-i)\ne0\big]
\le\frac1n\sum_{i=1}^n\big(1-e^{-\theta/i}+\tfrac1{n-i}\big).
\end{equation}
This expression vanishes as $n\rightarrow\infty$, hence $\tilde\alpha_{i^*(n)}(n)=0$ with high probability (this means that the part that differs between $\alpha$ and $\tilde\alpha$ does not appear elsewhere in $\alpha$).

Since $i^*(n)$ is the size of the final gap between 1s that appears in $(\xi_1,\ldots,\xi_n,1)$,
\begin{align}
\P\big[i^*(n)=k\big]
&=\P\big[\xi_n=0,\ldots,\xi_{n-k+2}=0,\xi_{n-k+1}=1\big]\\
&=\tfrac{n-1}{\theta+n-1}\tfrac{n-2}{\theta+n-2}\cdots\tfrac{n-k+1}{\theta+n-k+1}\tfrac\theta{\theta+n-k}.
\label{eqn:istar}
\end{align}
This vanishes for fixed $k$ as $n\rightarrow\infty$, so $i^*(n)\rightarrow\infty$ in probability. Since $e_i$ vanishes as $i\rightarrow\infty$, for any fixed $\epsilon>0$, $e_{i^*(n)}<\epsilon$ with high probability. Hence, with high probability,
\begin{equation}
\bigg|\frac{\eta_{n,n}}{\tilde \eta_{n,n}}-1\bigg|
=\bigg|\frac{\zeta_{i^*(n)}(\tilde\alpha_{i^*(n)}(n)+1)}{\zeta_{i^*(n)}(\tilde\alpha_{i^*(n)}(n))}-1\bigg|
=\bigg|\frac{\zeta_{i^*(n)}(1)}{\zeta_{i^*(n)}(0)}-1\bigg|
\le\frac{e_{i^*(n)}}\theta\le\frac\epsilon\theta.
\end{equation}
This proves that $\eta_{n,n}\stackrel p\sim\tilde \eta_{n,n}$.

For $\tilde \eta_{n,n}\stackrel p\sim\tilde \eta_{n,m}$, observe that with high probability, $\alpha_i(\infty)\le1$ for all $i>m$. Picking $m$ so that $e_i\le\frac\theta2$ for all $i>m$, and using the inequality $|\log(1-x)|\le2\log(1+x)$ for $x\ge\frac12$, with high probability,
\begin{equation}
\bigg|\log\frac{\tilde \eta_{n,n}}{\tilde \eta_{n,m}}\bigg|
\le\sum_{i>m}\big|\log\zeta_i(\tilde\alpha_i(n))\big|
\le2\sum_{i>m}\log\zeta^+_i(\alpha_i(\infty)).
\end{equation}
This is twice the tail of a series whose sum is $\log\eta^+_\infty$, which vanishes since $\eta^+_\infty<\infty$ almost surely, hence $\tilde\eta_{n,n}/\tilde\eta_{n,m}\stackrel p\rightarrow1$. For the same reason, we have $\eta_\infty/\eta_{\infty,m}\stackrel{as}\rightarrow1$. Finally, since $\P\big[\tilde\alpha_i(n)=\alpha_i(n)=\alpha_i(\infty)\big]\rightarrow1\vphantom{\stackrel p\rightarrow}$ for any fixed $i$, $\P\big[\tilde\eta_{n,m}=\eta_{n,m}=\eta_{\infty,m}\big]\rightarrow1$.

\vspace{6pt}\noindent\textbf{Step 3: Uniform integrability.} Since $\eta^{}_{n,m}\le\eta^+_{n,m}\le\eta^+_{n,n}$, it suffices to prove uniform integrability of $\eta^+_{n,n}$. Let $p_{n,i}=\P[i^*(n)=i]$. By Lemma \ref{lem:lcs} and independence of $\alpha_i(\infty)$,
\begin{align}
\E\big[\eta^+_{n,n}\big]
=\sum_{i=1}^np_{n,i}\E\big[\eta^+_{n,n}\big|i^*(n)=i\big]
&=\sum_{i=1}^np_{n,i}\E\bigg[\zeta_i^+(\tilde\alpha_i(n)+1)\prod_{j\ne i}\zeta_j^+(\tilde\alpha_j(n))\bigg]\\
&\le\sum_{i=1}^np_{n,i}\E\big[\zeta_i^+(\alpha_i(\infty)+1)\big]\prod_{j\ne i}\E\big[\zeta_j^+(\alpha_j(\infty))\big]\\
&=\E\big[\eta^+_\infty\big]\sum_{i=1}^np_{n,i}\frac{\E\big[\zeta^+_i(\alpha_i(\infty)+1)\big]}{\E\big[\zeta^+_i(\alpha_i(\infty))\big]}.
\end{align}
Let $C=\sum_\ell c_\ell<\infty$ and $D=\sum_\ell\ell c_\ell<\infty$. Since $\alpha_i(\infty)$ is Poisson with parameter $\theta/i$, $\E\big[\zeta^+_i(\alpha_i(\infty)+1)\big]=1+De_i/\theta$ and $\E\big[\zeta^+_i(\alpha_i(\infty))\big]=1+(\theta+Ce_i)/i\ge1$, hence
\begin{equation}
\sum_{i=1}^np_{n,i}\frac{\E\big[\zeta^+_i(\alpha_i(\infty)+1)\big]}{\E\big[\zeta^+_i(\alpha_i(\infty))\big]}
\le\sum_{i=1}^np_{n,i}\Big(1+\frac{De_i}\theta\Big)
=1+\frac D\theta\E\big[e_{i^*(n)}\big].
\label{eqn:etaratio}
\end{equation}
We previously showed that $e_{i^*(n)}$ is bounded and vanishes in probability, so (\ref{eqn:etaratio}) converges to 1, hence $\limsup\E[\eta^+_{n,n}]\le\E[\eta^+_\infty]$. But $\eta^+_{n,n}\ge\tilde\eta^+_{n,n}$, while $\tilde\eta^+_{n,n}\rightarrow\eta^+_\infty$ in $L^1$ by dominated convergence, thus
\begin{equation}
\limsup_{n\rightarrow\infty}\E\big[\big|\eta^+_{n,n}-\tilde\eta^+_{n.n}\big|\big]
=\limsup_{n\rightarrow\infty}\E[\eta^+_{n,n}]-\lim_{n\rightarrow\infty}\E[\tilde\eta^+_{n,n}]\le0.
\end{equation}
Hence, $\eta^+_{n,n}\rightarrow\eta^+_\infty$ in $L^1$ and is therefore uniformly integrable.
\end{proof}

\section{Conclusion}
We studied measures on partitions in terms of perturbations of the Ewens measure, and showed that under a certain condition, many limit theorems for the Ewens measure are universal. This unifies the proofs of the limit theorems under a single universality theorem, while simultaneously extending the class of measures for which they hold beyond the previous frontier of logarithmic combinatorial structures. An interesting direction for future work is to better understand the asymptotically Ewens condition, which will perhaps allow more cases to be unified under the universality theorem.




\setlength\bibsep{0cm}


\ACKNO{The author wishes to thank Persi Diaconis and Yunjiang Jiang for the many enlightening discussions on this topic.}


\end{document}